\documentclass[12pt]{article}
\usepackage[a4paper,left=2cm,right=1.2cm,top=3cm,bottom=4cm]
{geometry}
\usepackage{amsmath,amstext, amsthm, empheq, extpfeil,  
amsbsy,amssymb,marvosym,fancyhdr,graphicx,amscd,amsfonts,latexsym,delarray,stackrel,lineno,color,cite,appendix,xcolor,url,hyperref}
\usepackage{wasysym}
\usepackage{subfig}

\usepackage{srcltx}
\usepackage{mathrsfs}
\usepackage{float} 
\usepackage[all]{xy}
\usepackage{t1enc}
\usepackage{mathrsfs}
\usepackage{pifont}

\usepackage{mathpple}
\usepackage[T1]{fontenc}

\definecolor{Green}{rgb}{0,1,0}
\definecolor{Blue}{RGB}{0,0,191}
\definecolor{mathmodecolor}{RGB}{0,102,0}
\definecolor{keywordcolor}{RGB}{0,51,151}
\definecolor{sourcebackgroundcolor}{RGB}{255,247,223}
\definecolor{unixagred}{RGB}{255,0,0}
\definecolor{lightgray}{RGB}{191,191,191}
\definecolor{green}{RGB}{1,191,191}

\newcommand*\patchAmsMathEnvironmentForLineno[1]{%
  \expandafter\let\csname old#1\expandafter\endcsname\csname #1\endcsname
  \expandafter\let\csname oldend#1\expandafter\endcsname\csname end#1\endcsname
  \renewenvironment{#1}%
     {\linenomath\csname old#1\endcsname}%
     {\csname oldend#1\endcsname\endlinenomath}}%
\newcommand*\patchBothAmsMathEnvironmentsForLineno[1]{%
  \patchAmsMathEnvironmentForLineno{#1}%
  \patchAmsMathEnvironmentForLineno{#1*}}%
\AtBeginDocument{%
\patchBothAmsMathEnvironmentsForLineno{equation}%
\patchBothAmsMathEnvironmentsForLineno{align}%
\patchBothAmsMathEnvironmentsForLineno{flalign}%
\patchBothAmsMathEnvironmentsForLineno{alignat}%
\patchBothAmsMathEnvironmentsForLineno{gather}%
\patchBothAmsMathEnvironmentsForLineno{multline}%
}

\newtheorem{thm}{Theorem}[section]
\newtheorem{prop}[thm]{Proposition}

\newtheorem{lem}[thm]{Lemma}

\DeclarePairedDelimiter\floor{\lfloor}{\rfloor}

\def\A{{\mathbb A}}

\def\N{{\mathbb N}}

\def\Q{{\mathbb Q}}
\def\R{{\mathbb R}}
\def\Z{{\mathbb Z}}

\def\aarith{{\mathfrak A}}

\def\cO{{\mathcal O}}

\def\scal{{(\rnt,\cO)}}
\def\dar[#1]{\ar@<2pt>[#1]\ar@<-2pt>[#1]}
\def\qqq{\,,\,~\forall}

\newcommand{\ie}{{\it i.e.\/}\ }

\def\sin{{{\rm sin}}}
\def\cos{{{\rm cos}}}
\def\tan{{{\rm tan}}}

\def\bm2{{\rm Bmod^2}}
\def\b2{{\rm Bmod^{\mathfrak s}}}

\newcommand{\nil}[1]{}

\def\rnt{{[0,\infty)\rtimes{\N^{\times}}}}

\def\aarith{{\mathscr A}}
\def\scal{{(\rnt,\cO)}}
\def\scal1{{\hat \aarith}}
\def\scal2{{\mathscr S}}

\title
{Around Wilson's theorem}

\author{Alain Connes}

\begin{document}
\maketitle

\abstract{We study the series s(n,x) which is the sum for k from 1 to n of the square of the sine of the product x Gamma(k)/k, where x is a variable. By Wilson's theorem we show that the integer part of s(n,x) for x = Pi/2 is the number of primes less or equal to n and we get a similar formula for x a rational multiple of Pi. We show that for almost all x in the Lebesgue measure s(n,x) is equivalent to n/2 when n tends to infinity, while for almost all x in the Baire sense, 1/2 is a limit point of the ratio of s(n,x) to the number of primes less or equal to n. }
\section{Introduction}
Let  $\Pi(n)$ be the number of primes $p\leq n$.  A slight improvement on a formula\footnote{\begin{figure}[H]	\begin{center}
\includegraphics[scale=0.5]{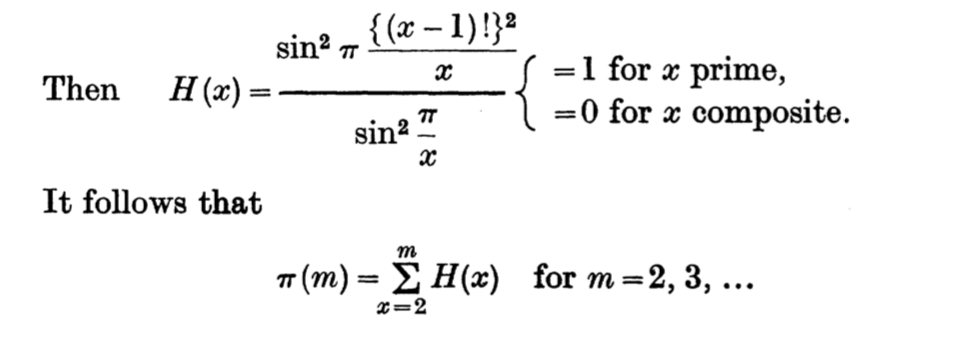}
\end{center}
\caption{The formula of Willans. \label{willans} }
\end{figure}} of Willans \cite{Willans} gives  a simple formula for $\Pi(n)$ as the integer part of the sum
\begin{equation}\label{firstf}
\sum _{k=1}^n \sin^2\left(\frac{\pi  \Gamma (k)}{2k}\right)
\end{equation}
When one tries to compute naively the right hand side one finds that it requires an increasing precision on the numerical value of the number $\pi$ whose first $2500$ decimals are needed to compute $\Pi(n)$ for $n$ of the order of a thousand. F. Villegas suggested to replace $\pi$ by a variable and analyse the dependence on $x$ in the above sequence. 
Thus for $n>1$ an integer and $x\in \R$,  let 
	\begin{equation}\label{functionf}
		s(n,x):=\sum _{k=1}^n \sin^2\left(\frac{x  \Gamma (k)}{k}\right)
	\end{equation}
	We shall show below that the dependence on $x\in \R$ is quite interesting inasmuch as, due to the lacunary nature of the sequence $\frac{\Gamma (k)}{k}$, the terms of the sum \eqref{firstf} are essentially independent random variables when suitably understood as functions on an almost periodic compactification $G$ of $\R$. This  gives, by the proof of the strong law of large numbers, that for almost all $x\in \R$ in the sense of the Lebesgue measure one has when $n\rightarrow\infty$ that $s(n,x)\sim \frac n2$. The interesting fact is that for the other natural notion of ``generic" real number, namely the one provided by the Baire theory of dense countable intersections of open sets, it is a totally different behavior of the sequence $s(n,x)$ which is generic:  we show in Theorem \ref{generic}
that for generic $x\in \R$,  the quotients $\frac{s(n,x)}{\Pi(n)}$  get arbitrarily close to $\frac 12$, \ie $\frac 12$ is a limit point of the sequence
$$
\frac 12 \in \lim_{n\to \infty} \frac{s(n,x)}{\Pi(n)}.
$$
Generically this sequence will also have $\infty$ as a limit point and will oscillate wildly. But for rational multiples of $\pi$ the sequence $s(n,x)$ behaves like the product of $\Pi(n)$ by the  rational number\footnote{$\mu$ is the M\" oebius function and $\phi$ the Euler totient function.} $\frac{1}{2} -\frac{\mu(b)}{2 \phi(b)}$ which only depends upon the denominator $b>1$ of the irreducible fraction $x= \frac ab \pi$ as a multiple of $\pi$ (see Proposition \ref{behavepi}). 

\section{$\Pi(n)$ and sum of squared sines}

\begin{figure}[H]	\begin{center}
\includegraphics[scale=0.7]{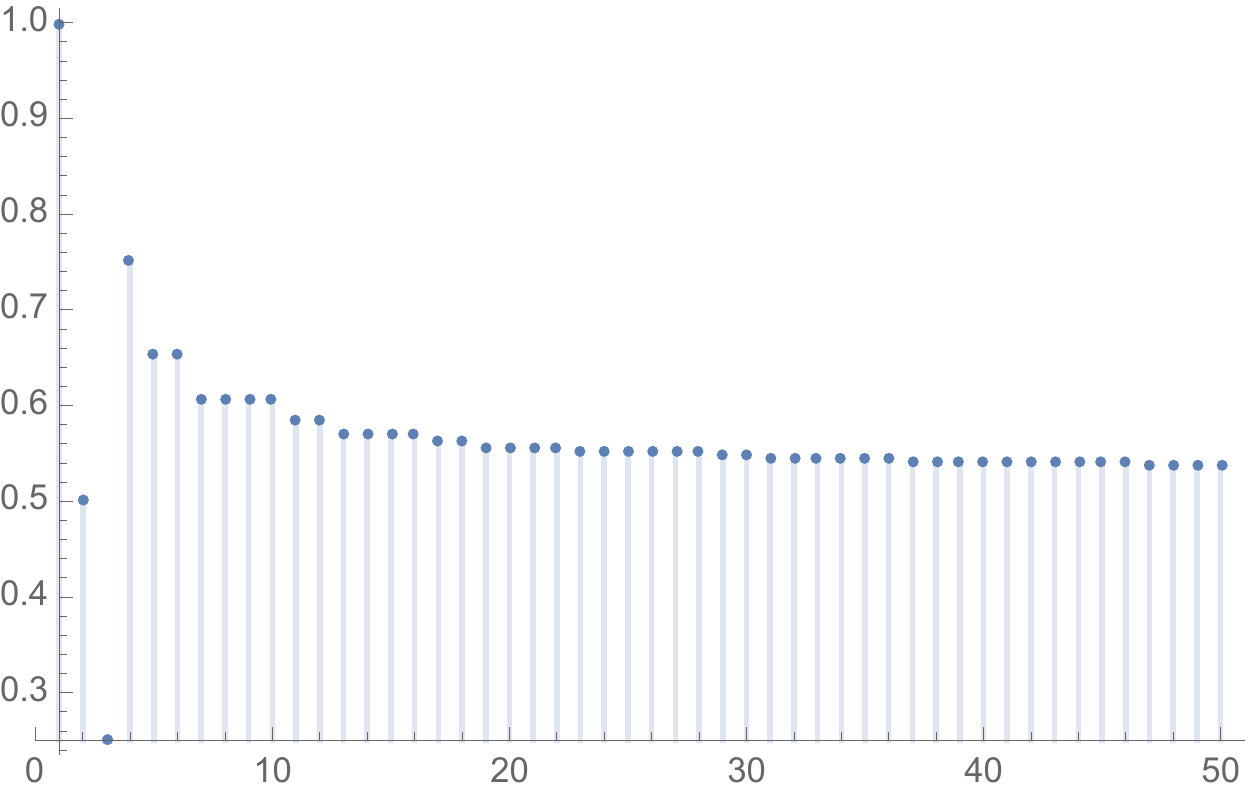}
\end{center}
\caption{Graph of $s(n,\frac \pi 2)-\Pi(n)$ for $1<n\leq 50$. \label{feyn2} }
\end{figure}

We start with the following variant of the formulas of Willans \cite{Willans}.
\begin{prop}\label{formpi}
	Let $n>1$ be an integer then $\Pi(n)$ is the integer part of $s(n,\frac \pi 2)$. 
\end{prop}
\proof 
For $k>4$ not prime the quotient $(k-1)!/k$ is an even integer. Thus in that case one has $$\sin^2 \left(\frac{\pi  \Gamma (k)}{2k}\right)=0$$
For $k=p>2$ prime, the residue of $(p-1)!$ modulo $2p$ is $p-1$ by Wilson's theorem and the evenness of $p-1$. Thus for $p>2$ prime, 
$$
\sin^2 \left(\frac{\pi  \Gamma (p)}{2p}\right)=\sin^2 \left(\frac{\pi  (p-1)}{2p}\right)=\cos^2 \left(\frac{\pi}{2p}\right).
$$
One has 
$$
1\geq \cos^2(x)\geq 1- x^2\qqq x\in \R
$$
It follows that $\delta(n)=s(n,\frac \pi 2)-\Pi(n)$ is a decreasing function of $n>4$ and that for $m>n$, $$\delta(n)-\delta(m)=\sum_{p\,\text{ prime}, \, n< p\leq m}\left(1-\cos^2 \left(\frac{\pi}{2p}\right)\right)\leq \frac{\pi^2}{4}\sum_{n< p\leq m} \frac{1}{p^2}$$
The series $\sum \frac{1}{p^2}$ is convergent and one has the bound $
\frac{\pi^2}{4}\sum_{p> 50} \frac{1}{p^2}<0.0498448
$ for the (larger) sum over integers, 
while $s(50,\frac \pi 2)-\Pi(50)\sim 0.539005$. It follows that $$s(n,\frac \pi 2)\in [\Pi(n)+0.5-0.05),\Pi(n)+0.6]\subset [\Pi(n),\Pi(n)+1)$$ for any $n>50$, and (see Figure \ref{feyn2}) $s(n,\frac \pi 2)\in [\Pi(n),\Pi(n)+1)$ for any $n>1$ which gives the required result.\endproof

The general term $\sin^2\left(\frac{x  \Gamma (k)}{k}\right)$ of \eqref{functionf} depends on the knowledge of $x$ up to an $\epsilon$ of the order of 
$$
dx\simeq \left(\frac k e\right)^{-k+2}
$$
Thus for instance to get the required precision around $k=945$ one needs the first $2400$ decimals of $\pi$.
\begin{figure}[H]	\begin{center}
\includegraphics[scale=0.4]{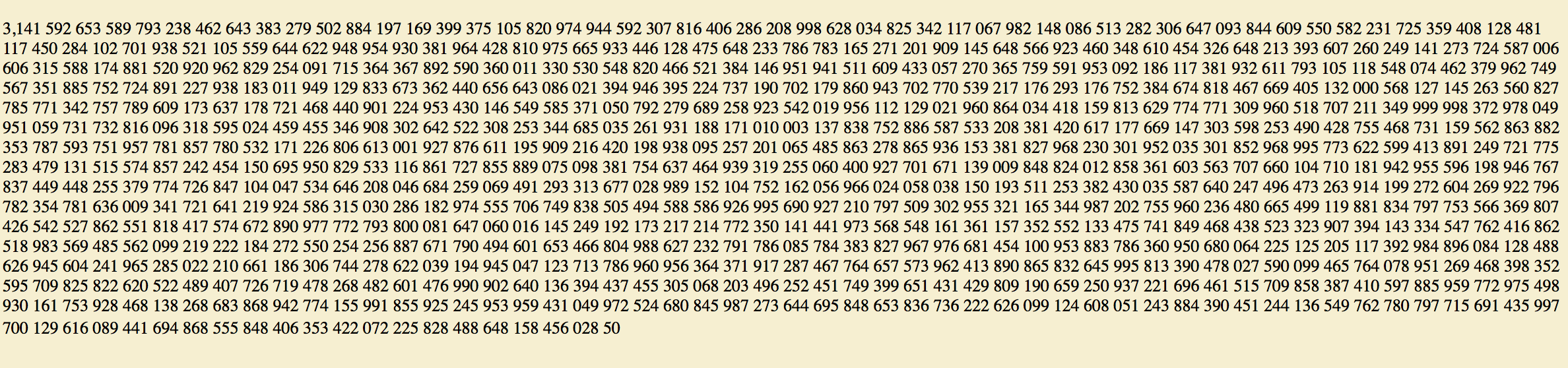}
\end{center}
\caption{$2400$ decimals of $\pi$. \label{pi} }
\end{figure}
\section{Rational multiples of $\pi$} 
We investigate  the behavior of the sequence $s(n,x)$ when $x$ is a rational multiple of $\pi$. We let $x:=\frac{a \pi}{ b}$ where $a$ and $b$ are relatively prime. For $k$ large enough the term $\Gamma(k)/k$ is divisible by $b$ in the sense that the $p$-adic valuation of $\Gamma(k)/k$ is larger than that of $b$ for all prime divisors of $b$. This can be seen using the Legendre formula for the $p$-adic valuation which gives
$$
v_p(\Gamma(k)/k)\geq\sum_{\ell\geq 1} \floor*{(k-1)/p^\ell}-\floor*{\log_p(k)}
$$
 Since the prime factors of $b$ are fixed there exists $k_0<\infty$ such that $\Gamma(k)/k$ is divisible by $b$ (in the above sense) for all $k>k_0$. Thus if $k>k_0$ is not prime the product $\frac{x  \Gamma (k)}{k}$ is an integer multiple of $\pi$ and $\sin^2\left(\frac{x  \Gamma (k)}{k}\right)=0$. When $k>k_0$ is a prime, the integer $c=\frac{  \Gamma (k)}{b}$ is such that $bc=-1$ modulo $k$ by Wilson's theorem. One gets in this case 
$$
\sin^2\left(\frac{x  \Gamma (k)}{k}\right)=\sin^2\left(\frac{\pi a \Gamma (k)}{bk}\right)=\sin^2\left(\frac{\pi  a c}{k}\right)
$$
and the right hand side only depends upon the residue of $a$ and of $c$ modulo $k$. Since $b$ and $k$ are relatively prime there exists $u\in \{1, \ldots b-1\}$ which is the inverse of $k$ modulo $b$. Thus let $m\in \N$ such that $ku-1=bm$. One then has $bm=-1$ modulo $k$ and it follows that $m=c$ modulo $k$.
This gives 
$$
\sin^2\left(\frac{\pi  ac}{k}\right)=\sin^2\left(\frac{\pi a  m}{k}\right)=\sin^2\left(\frac{\pi a(ku-1)}{bk}\right)
=\sin^2\left(\frac{\pi a u}{b}\right)+\sin \left(\frac{a \pi }{b k}\right) \sin \left(\frac{a \pi }{b k}-\frac{2 a u \pi }{b}\right)
$$
Thus, since $$\vert\sin \left(\frac{a \pi }{b k}\right) \sin \left(\frac{a \pi }{b k}-\frac{2 a u \pi }{b}\right)\vert\leq \vert\sin \left(\frac{a \pi }{b k}\right)\vert =0(1/k) $$ the asymptotic behavior of $s(n,x)$ only depends, up to a term of the order of $\log\log n$ due to the $\sum 1/k$ over primes less than $n$, upon the residues of the primes modulo $b$. Thus by Dirichlet's theorem, in the strong form due to La Vallée Poussin (see \cite{Prachar}, V \S 7), one gets 
\begin{prop}\label{behavepi}
	Let $x$ be a rational multiple of $\pi$, $x=\frac{\pi a}{b}$ with $a,b$ relatively prime, $b>1$, then 
	\begin{equation}\label{ratcase}
		s(n,x)\sim \frac{\Pi(n)}{\phi(b)}\sum_{v\in (\Z/b\Z)^*}\sin^2\left(\frac{\pi v}{b}\right)
=\Pi(n)\left(\frac{1}{2} -\frac{\mu(b)}{2 \phi(b)}\right)
	\end{equation}
\end{prop}
\proof It remains to show the second equality in \eqref{ratcase}. It follows from $ \sin^2 (v)=\frac{1-\cos(2 v)}{2}$ and the fact that the sum of primitive roots of unity of order $b$ is the M\"obius function $\mu(b)$.\endproof 

\section{Generic behavior of $s(n,x)$}
This suggests to investigate the general behavior of the sequence \eqref{functionf}. It is given by the following result which shows that almost everywhere in measure theory (for the Lebesgue measure) the behavior is Gaussian and $s(n,x)\simeq \frac n2$. But generically at the topological level which means on a dense countable intersection of open sets, the sum grows far more slowly and behaves like $\frac{\Pi(n)}{2}$ in the  weak sense that $\frac 12$ is a limit point of the sequence $\frac{s(n,x)}{\Pi(n)}$. In fact it oscillates wildly since $\infty$ is also a limit point of that sequence. Note that the two behaviors are exclusive of each other but this is not a contradiction.
\begin{thm}\label{generic}
$(i)$~For almost all $x\in \R$ the sequence 	$s(n,x)$ of \eqref{functionf} has the Gaussian behavior
$$
s(n,x)\simeq \frac n2 
$$
$(ii)$~For generic $x\in \R$ (\ie on a dense countable intersection of open sets) one has 
$$
\frac 12 \in \lim_{n\to \infty} \frac{s(n,x)}{\Pi(n)}
$$
\end{thm}
\proof $(i)$~Let $G$ be the compact group  projective limit of the compact groups $G_n:=\R/n\Z$ under the natural morphisms 
\begin{equation}\label{projgn}
  \gamma_{n,m}:G_m\to G_n, \qquad \gamma_{n,m}(x+m\Z)=x+n\Z,\qquad\forall n\vert m.
  \end{equation}
  One has a natural isomorphism, with $\A_\Q$ the adeles of the global field $\Q$,
 $$
 G=\varprojlim (\Z/n\Z\times \R)/\Z =  (\hat\Z\times \R)/\Z= \A_\Q/(\Q,+).
  $$
 The Pontrjagin dual of $G$ is identified with the discrete additive group $\Q$ of rational numbers by associating to $r\in \Q$ the character $\alpha_r$ of $G$ specified by its restriction to the dense subgroup $\R$, range of the  homomorphism $\R\ni t\mapsto a(t)=(0,t)\in \A_\Q/(\Q,+)$ of adeles with $0$ non-archimedean component
 $$
 \alpha_r(a(t)):=e^{2 \pi i r t}.
 $$
Next, using $1-\cos(2 x)=2\, \sin^2 (x)$ we get with $r=\frac{ \Gamma (k)}{k}\in \Q$ the equality 
\begin{equation}\label{charac}
  \sin^2\left(\frac{x  \Gamma (k)}{k}\right)=
  \frac 12-\frac 14  \alpha_r\left(\frac x \pi\right) 
  -\frac 14  \alpha_{-r}\left(\frac x \pi\right).
  \end{equation}
Thus with the basic function defined on $G$ by $$X(x):=
-\frac 14  \alpha_1(x) 
  -\frac 14  \alpha_{-1}(x)$$
we get that the general term of the sum $s(n,x)$ is simply $\frac 12+X_k\left(\frac x \pi\right)$ where
$$
X_k(x):=X(r(k)x),\ \ r(k):=\frac{ \Gamma (k)}{k} \in \Q^\times
$$
The multiplication by elements of $\Q^\times$ defines automorphisms of $G$. One has the orthogonality relation of characters which implies since the rationals $\frac{\Gamma (k)}{k}$ are distinct, they form (for $k>1$) a strictly increasing sequence (the first ones are $\left\{1,\frac{1}{2},\frac{2}{3},\frac{3}{2},\frac{24}{5}\right\}$) that the random variables $X_k$ on the probability space $G$ equipped with its normalized Haar measure, are equidistributed and essentially independent inasmuch as,
$
\int_G X_k(x) X_\ell (x)dx =0 \qqq k\neq \ell 
$
and that one controls the $4$-th moment as follows 
\begin{equation}\label{moments}
\int_G\vert \sum_1^n   X_k(x) \vert^4 dx\leq C n^2
\end{equation}
since $\vert \sum_1^n   X_k(x) \vert^4= \sum_{1\leq k_j\leq n}   X_{k_1}(x)X_{k_2}(x)X_{k_3}(x)X_{k_4}(x) $ and the number of solutions of the equation 
\begin{equation}\label{cancell}
\pm r(k_1)\pm r(k_2)\pm r(k_3)\pm r(k_4)=0, \ \ k_j\in \{1,\ldots ,n\}
\end{equation}
is of the order of $n^2$ due to the lacunary nature (\!\!\cite{SZ}) of the sequence $r(k)$. Indeed for $k>4$ one has $r(k+1)>3r(k)$ and  \eqref{cancell} is possible only if the largest $k_j$ appears at list twice (and with opposite signs) and the remaining $k_j$ are equal, which gives $n^2$ as a bound on the number of solutions.

Thus one has, for any $\epsilon >0$ that 
$$
\int_G\vert \frac 1n\sum_1^n X_k(x)\vert^4\leq C/n^2, \ \ \vert\{x\in G\mid \vert \frac 1n\sum_1^n X_k(x)\vert> \epsilon\}\vert\leq C/n^2 \epsilon^{-4}
$$
and the Borel-Cantelli lemma  applies and shows that the subset $E\subset G$ defined by 
$$
E:=\{x\in G\mid \frac 1n\sum_1^n X_k(x)\to 0\}
$$
is of measure $1$. Since $\R\subset G$ is of measure $0$ we cannot yet get $(i)$ but it will follow from the invariance of $E$ under the translation by the subgroup $\hat\Z\subset G$. To see this we use the equality for $k$ non-prime
\begin{equation}\label{charac1}
  X_k(x+u)=X_k(x) \qqq u\in \hat\Z
  \end{equation}
  which follows from the integrality of $\frac{ \Gamma (k)}{k}$ and the periodicity of the cosine: 
$$
\cos \frac{ 2\pi (x+1)\Gamma (k)}{k}=\cos \frac{ 2\pi x\Gamma (k)}{k}$$
Thus one gets the  same equality for the closure $\hat\Z\subset G$. It follows that
$$
\vert \frac 1n \sum_1^n X_k(x+u)- \frac 1n \sum_1^n X_k(x)\vert \leq \frac{\Pi(n)}{n}\qqq u\in \hat\Z
$$
and this suffices to show that $E$ is invariant under the translation by the subgroup $\hat\Z\subset G$. The image of $x\in G$ in the quotient $G/\hat\Z=\R /\Z$ thus suffices to decide if $x\in E$ and it follows that almost all elements of $\R\subset G$ are in $E$. Finally the Gaussian behavior follows from the results of \cite{SZ} on lacunary trigonometric series. \newline
$(ii)$~When $x\in \pi \Q$ is a rational multiple of $\pi$ one applies Proposition \ref{behavepi}. For $b\to \infty$ one has the equidistribution 
$$
\frac{1}{\phi(b)}\sum_{v\in (\Z/b\Z)^*}\sin^2\left(\frac{\pi v}{b}\right)=\frac{1}{2} -\frac{\mu(b)}{2 \phi(b)}\to \frac 12.
$$
This shows that, for $\epsilon >0$, the following countable intersection of  open sets is dense in $\R$
$$
W(\epsilon):=\cap_m\cup_{n\geq m}\{x\in \R\mid s(n,x)\in((1-\epsilon) \Pi(n)/2,(1+\epsilon) \Pi(n)/2)\}
$$
and for $x\in W(\epsilon)$ one has $\lim_{n\to \infty} \frac{2 s(n,x)}{\Pi(n)}\cap [1-\epsilon,1+\epsilon]\neq \emptyset$ which gives the required conclusion after intersecting the $W(\epsilon)$ for $\epsilon=\frac 1a$, $a\to \infty$ which still gives  a dense countable intersection of  open sets by Baire's Theorem \cite{baire}. \qed


\begin{thebibliography}{99}
 
 \bibitem{baire} R. Baire. {\em Sur les fonctions de variables réelles}. Ann. di Mat., 3:1--123, (1899).
 
 

\bibitem{Prachar} K.~Prachar, {\em Primzahlverteilung}. (German) Springer-Verlag, Berlin-Göttingen-Heidelberg, (1957). x+415 pp.

\bibitem{SZ} 
R.~Salem, A.~Zygmund, {\em On Lacunary Trigonometric Series}.
Proceedings of the National Academy of Sciences of the United States of America
Vol. 33, No. 11 (Nov. 15, 1947), pp. 333--338

\bibitem{Willans} C.~P.~Willans, {\em On formulae for the nth prime number}. Math. Gaz. 48 (1964) 413--415.

\end{thebibliography}
\end{document}